\documentclass[10pt,reqno]{amsart}
\usepackage{amssymb}
\usepackage[all]{xy}

\theoremstyle{plain}
\newtheorem{prop}{Proposition}[section]
\newtheorem{theo}[prop]{Theorem}
\newtheorem{coro}[prop]{Corollary}
\newtheorem{lemm}[prop]{Lemma}
\theoremstyle{remark}
\newtheorem{rema}[prop]{Remark}
\theoremstyle{definition}

\newtheorem*{exam}{Example}
\newtheorem{assu}[prop]{Assumption}

\numberwithin{equation}{section}

\newcommand{\eqto}{\stackrel{\lower1.5pt\hbox{$\scriptstyle\sim\,$}}\to}

\newcommand{\clap}[1]{\hbox to 0pt{\hss #1\hss}}

\def\cA{{\mathcal A}}

\def\cE{{\mathcal E}}

\def\cO{{\mathcal O}}
\def\cL{{\mathcal L}}

\def\Gal{{\rm Gal}}
\def\inv{{\rm inv}}
\def\disc{{\rm disc}}

\def\dv{{\rm div}}

\def\A{{\mathbb A}}
\def\C{{\mathbb C}}

\def\PP{{\mathbb P}}
\def\Q{{\mathbb Q}}
\def\G{{\mathbb G}}
\def\R{{\mathbb R}}

\def\Z{{\mathbb Z}}
\def\C{{\mathbb C}}

\def\Pic{{\rm Pic}}

\def\Hom{{\rm Hom}}
\def\Br{{\rm Br}}
\def\im{{\rm im}}

\begin{document}
\title[Brauer--Manin obstructions]
{Effectivity of Brauer--Manin obstructions}
\author{Andrew Kresch}
\address{
  Institut f\"ur Mathematik,
  Universit\"at Z\"urich,
  Winterthurerstrasse 190,
  CH-8057 Z\"urich, Switzerland
}
\email{andrew.kresch@math.unizh.ch}
\author{Yuri Tschinkel}
\address{
  Mathematisches Institut,
  Georg-August-Universit\"at G\"ottingen,
  Bunsenstrasse 3-5,
  D-37073 G\"ottingen, Germany
}
\email{yuri@uni-math.gwdg.de}

\date{December 21, 2006}
\thanks{The first author was supported by
an Alexander von Humboldt Foundation Research Fellowship.
The second author was supported by the NSF}

\begin{abstract}
We study Brauer--Manin obstructions to the Hasse principle and to
weak approximation, with special regard to effectivity questions.
\end{abstract}
\maketitle

\section{Introduction}
\label{sec:introduction}

Let $k$ be a number field, $X$ be a smooth projective
(geometrically irreducible) variety
over $k$, and $X(k)$ its set of $k$-rational points.
An important problem in arithmetic geometry is to find an effective procedure
to determine whether $X(k)\ne \emptyset$.
A necessary condition is that $X(k_v)\ne \emptyset$ for all completions
$k_v$ of $k$.
This condition can be tested effectively and easily, given the defining equations
for $X$.
One says that $X$ \emph{satisfies the Hasse principle} when
\begin{equation}
\label{eqn:hasseprinciple}
X(k)\ne \emptyset \Leftrightarrow X(k_v)\ne \emptyset \,\,\,\forall v.
\end{equation}

When $X$ is a quadric hypersurface (of arbitrary dimension) over
the rational numbers, the
validity of \eqref{eqn:hasseprinciple} was established in 1921 by
Hasse in his doctoral thesis \cite{hasse1}.
The statement \eqref{eqn:hasseprinciple} was proposed as a
\emph{principle} by Hasse in 1924 \cite{hasse3}, where it was
proved to hold for quadric hypersurfaces
over arbitrary number fields.
Hasse's main insight was to relate the existence of solutions to equations
over a number field to existence of solutions over its completions,
i.e., the $v$-adic numbers, which had been
introduced and developed into a theory by his
thesis advisor Hensel \cite{hensel08}.
In fact, Hensel had studied $v$-adic solutions to quadratic equations
(see \cite{hensel13} Chapter 12),
obtaining necessary and sufficient conditions for local solvability.
Earlier, Minkowski had defined a complete system of invariants of
quadratic forms over the rationals, one at each prime $p$ \cite{minkowski};
the Hasse principle for quadratic forms \cite{hasse2} (or
in other settings)
is often referred to as the Hasse--Minkowski principle.

A related problem in arithmetic geometry
is to find $k$-rational points on $X$ matching local data,
i.e., determining whether or not $X(k)$ is dense in the adelic space
$$X(\A_k)=\prod_v X(k_v).$$
(The adelic space is equipped with the product topology.)
In this case one says that $X$ \emph{satisfies weak approximation}.

The Hasse principle, and weak approximation, are known to fail for general
projective varieties, e.g., cubic curves and cubic surfaces.
Counterexamples to the Hasse principle appeared as early as
1880 \cite{pepin}. (For a discussion, and proofs of the claims that
appeared at that time, see \cite{lemmermeyer}.)
By the early 1940's it was well established that genus $1$ curves may fail
to satisfy the Hasse principle \cite{lind}, \cite{reichardt}.

All known obstructions to the Hasse principle and weak approximation
are based on the \emph{Brauer--Manin obstruction} defined by Manin
\cite{manin}
or reduce to this after finite \'etale covers of $X$ \cite{skoro}.
The Brauer--Manin obstruction is based on the Brauer group
$\Br(X)$ of $X$
and class field theory for $k$.
It cuts out a subspace $X(\A_k)^{\Br}$ of the adelic space
$X(\A_k)$ with the property
that $X(k)\subset X(\A_k)^{\Br}$.
In particular, if $X(\A_k)\ne\emptyset$ and $X(\A_k)^{\Br}=\emptyset$, then
$X(k)=\emptyset$ and $X$ fails to satisfy the Hasse principle.
If $X(k)\ne\emptyset$ and
$X(\A_k)^{\Br}\ne X(\A_k)$ then weak approximation fails for $X$.
In the former case, we speak of a \emph{Brauer--Manin obstruction
to the Hasse principle}, and in the latter case, of a \emph{Brauer--Manin
obstruction to weak approximation}.
We explain this in Section \ref{sec:Brauer}.

For geometrically rational surfaces, one expects that
the Brauer--Manin obstruction is the only obstruction to the Hasse principle
and to weak approximation \cite{cts}.
However, the explicit computation of this obstruction is a nontrivial
task, even in such concrete classical cases as cubic surfaces over $\Q$
\cite{ctks}, \cite{kt}.

In this note, we prove (Theorem \ref{iscomputable}) that
there is a general procedure for computing this
obstruction, provided the geometric Picard group is finitely generated,
torsion free, and known explicitly by means of cycle representatives
with an explicit Galois action.
See Section \ref{sec:prelim} for a precise description of the
required input data.
In particular this procedure is applicable to all del Pezzo surfaces
(see Section \ref{sec:prelim}).

The procedure is presented in Section \ref{sec:procedure}.
It involves several steps, which we now summarize.
A splitting field (a finite extension of $k$
over which all the cycle representatives are defined) must be chosen;
this is taken to be Galois.
The first Galois cohomology group of the geometric Picard group can be
identified abstractly with a Brauer group.
Generators, which are $1$-cocycles for group cohomology, are computed.
The next step is to obtain cocycle data for Brauer group elements
from these generators, i.e., $2$-cocycles of rational functions on $X$.
With this, the computation of the subspace of $X(\A_k)$ cut out by
these Brauer group elements is carried out.
Some results in this direction have been discussed in
\cite{bsd}, \cite{bright}.
A \texttt{Magma} package for degree 4 del Pezzo surfaces is available
\cite{logan}.

\section{Preliminaries}
\label{sec:prelim}
In this section we introduce notation and give details of the
input data required for the algorithm of Section \ref{sec:procedure}.
We fix a number field $k$ and a smooth projective geometrically irreducible
variety $X$ over $k$.
When $K/k$ is a field extension, we write $X_K$ for the base change of
$X$ to $K$.
We write $\Gal(K/k)$ for the Galois group when the extension is normal.
The Picard group of $X$ is denoted $\Pic(X)$.

\begin{assu}
\label{assumption1}
We suppose that we are given explicit equations defining
$X$ in $\PP^N$, i.e.,
generators $f_1$, $\ldots$, $f_r$, of the ideal
$\mathcal{J}=\mathcal{J}(X)\subset k[x_0,\ldots,x_N]$.
We assume that $X(\A_k)\ne\emptyset$, that $\Pic(X_{\bar k})$ is
torsion free, and
that the following are specified:
\begin{enumerate}
\item a collection of codimension one geometric cycles
$D_1$, $\ldots$, $D_m\in Z^1(X_{\bar k})$
whose classes generate $\Pic(X_{\bar k})$,
i.e., there is an exact sequence of abelian groups
\begin{equation}
\label{picseq}
0\to R\to \bigoplus_{i=1}^m \Z\cdot [D_i]\to
\Pic(X_{\bar k})\to 0;
\end{equation}
\item the subgroup of relations $R$;
\item a finite Galois extension $K$ of $k$,
over which the $D_i$ are defined, with known Galois group
$$G:=\Gal(K/k);$$
\item the action of $G$ on $\Pic(X_{\bar k})$.
\end{enumerate}
For simplicity, we assume that the cycles in (1) are effective,
and the collection of cycles is closed under the Galois action.
\end{assu}

We adopt the convention that the Galois action on the splitting field is a left
action, written $\alpha\mapsto {}^g\!\alpha$ ($g\in G$, $\alpha\in K$).
Hence there is an induced right action of $G$ on $X_K$, which we denote
by $a_g\colon X_K\to X_K$.
For a divisor $D$ on $X_K$ we denote $a_g^*D$ by ${}^g\!D$.
The action of $G$ on $\Pic(X_{\bar k})$, mentioned in (4), is the
action of pullback by $a_g$, meaning that $g\in G$
sends the class of $D$ to the class of ${}^g\!D$.

\begin{exam}
Let $X$ be a del Pezzo surface of degree $d\le 4$.
(Note that 
Hasse principle and weak approximation hold for $d\ge 5$ \cite{maninbook}.)
It is known that $X_{\bar k}$ is isomorphic to a
blow-up of $\PP^2$ in $9-d$ points.
For $d=3$, $4$, the anticanonical class $-\omega_X$ gives an embedding
$X\hookrightarrow \PP^d$; this supplies the ideal $\mathcal{J}$.
When $d=2$ we get an embedding
$X\hookrightarrow \PP(1,1,1,2)\hookrightarrow \PP^6$
from $-2\omega_X$, and when $d=1$ we get a embedding $X\hookrightarrow
\PP(1,1,2,3)\hookrightarrow \PP^{22}$ by $-6\omega_X$.

There is a finite collection of exceptional curves on $X_{\bar k}$,
with explicit equations.
For instance, when $d=3$, $X\subset \PP^3$
is a (smooth) cubic surface, and \cite{sousley} gives a procedure
for the determination of the lines on $X$.

The classes of the exceptional curves
generate $\Pic(X_{\bar k})$.
The number $n_d$ of these curves is given in the following table:
$$
\begin{array}{c|cccc}
d  &  1&2&3&4\\
\hline
n_d &240&56&27&16
\end{array}
$$
We may take $m=n_d$ in the exact sequence \eqref{picseq}.
Furthermore, we know that $\Pic(X_{\bar k})$ is isomorphic to
$\Z^{10-d}$.
Intersection numbers of curves on $X$ are readily computed,
and $R$ can be obtained using the fact that the intersection pairing
on $\Pic(X_{\bar k})$ is nondegenerate.
The Galois group $G=\Gal(K/k)$ acts by permutations on the set of
exceptional curves, and thus on their classes in
$\Pic(X_{\bar k})$.
\end{exam}

In this paper, we make use of computations of group cohomology,
cf.\ \cite{brown}.
For $G$ a finite group,
any resolution of $\Z$ as a $\Z[G]$-module (by finite free
$\Z[G]$-modules) defines the cohomology $H^i(G,M)$
of a $G$-module $M$,
as the $i$th cohomology of the complex obtained by
applying the functor $\Hom(-,M)$.
We write $(\mathsf{C}^\bullet(M),\partial)$ for this complex.
(This complex depends on the choice of resolution, which we suppress
in the notation; the cohomology is an invariant of $M$.)
We write $M^G=H^0(G,M)$ for the submodule of $M$ of
$G$-invariant elements.

\section{The Brauer--Manin obstruction}
\label{sec:Brauer}
In this section we recall basic facts about the Brauer groups of
fields and schemes.
The Brauer group $\Br(F)$ of a field $F$ is the group of
equivalence classes of central simple algebras over $F$.
The Brauer group $\Br(X)$ of a Noetherian scheme $X$,
defined by Grothendieck \cite{gb}, is the
group of equivalence classes of
sheaves of Azumaya algebras on $X$.
It injects naturally into the
cohomological Brauer group $\Br'(X)$, defined
as the torsion subgroup of $H^2(X,\G_m)$,
\'etale cohomology of the sheaf $\G_m$
of invertible regular functions
(multiplicative group scheme).
For $X$ projective over a field
(or more generally, arbitrary
$X$ possessing an ample invertible sheaf),
$\Br(X)=\Br'(X)$ by a result of Gabber
(see \cite{dejong}).

By class field theory (see, e.g., \cite{cf}),
the Brauer group of a number field $k$ fits in an exact sequence
\begin{equation}
\label{Brkseq}
0\to \Br(k)\longrightarrow \bigoplus_v \Br(k_v) \stackrel{\inv}\longrightarrow
\Q/\Z\to 0
\end{equation}
where the direct sum is over completions $v$ of $k$.
The Brauer groups of the fields $k_v$ are known by local class field theory.
More precisely, there is a \emph{local invariant}
\begin{equation}
\label{eq:localinvariants}
\inv_v\colon \Br(k_v)\eqto
\left\{
\begin{array}{ccl}
\Q/\Z &\text{when}& v\nmid \infty\\
(\frac{1}{2}\Z)/\Z & \text{when} & k_v=\R\\
0 & \text{when} & k_v=\C
\end{array}
\right.
\end{equation}
In \eqref{Brkseq},
$$\inv=\sum_v \inv_v.$$
Define
\begin{equation}
\label{defXBr}
X(\A_k)^\Br = \bigl\{\,(x_v)\in X(\A_k)\,\big|\,
\sum_v \inv_v(A(x_v))=0\,\,\,\forall A\in \Br(X)\,\bigr\}.
\end{equation}
By \eqref{Brkseq},
$X(k)\subset X(\A_k)^\Br$.

The Leray spectral sequence
\begin{equation}
\label{sseq}
H^p(\Gal(\bar k/k),H^q(X_{\bar k},\G_m))\Rightarrow H^{p+q}(X,\G_m)
\end{equation}
gives rise to an exact sequence
\begin{align}
\label{fromleray}
\begin{split}
0\to {}&\Pic(X)\to \Pic(X_{\bar k})^{\Gal(\bar k/k)}\to{}\\
&\quad\Br(k)\to\ker(\Br(X)\to\Br(K_{\bar k}))\to
H^1(\Gal(\bar k/k),\Pic(X_{\bar k}))\to 0.
\end{split}
\end{align}
Our assumptions imply (see Remark \ref{remaiso} below)
\begin{equation}
\label{BralgisH1}
\ker(\Br(X)\to \Br(X_{\bar k}))/\Br(k)\cong H^1(\Gal(\bar k/k),\Pic(X_{\bar k})).
\end{equation}
This kernel $\ker(\Br(X)\to \Br(X_{\bar k}))$ is known as the
\emph{algebraic part} of the Brauer group.
We define
\begin{equation}
\label{defXBralg}
X(\A_k)^{\mathrm{Br.alg}} = \bigl\{\,(x_v)\in X(\A_k)\,\big|\,
\sum_v \inv_v(A(x_v))=0\,\,\,\forall A\in \ker(\Br(X)\to \Br(X_{\bar k}))
\,\bigr\}.
\end{equation}

\begin{rema}
\label{remacosets}
By virtue of the sequence \eqref{Brkseq} it suffices in
\eqref{defXBr} and \eqref{defXBralg} to consider
one representative from each $\Br(k)$-coset in $\Br(X)$.
\end{rema}

\begin{rema}
\label{remaiso}
Since $X(\A_k)\ne \emptyset$,
\begin{itemize}
\item[(i)] $\Br(k)\hookrightarrow \Br(X)$;
\item[(ii)] $\Pic(X)\cong \Pic(X_{\bar k})^{\Gal(\bar k/k)}$.
\end{itemize}
Indeed, each of the homomorphisms
$\Br(k_v)\to \Br(X_{k_v})$ is split because points exist locally, and
fact (i) follows from the exact sequence \eqref{Brkseq}.
This implies the vanishing of the edge homomorphism
in \eqref{fromleray}, which in turn implies (ii).
\end{rema}

By assumption, the field $K$ is chosen so that
\begin{equation}
\label{PicXKPicXkbar}
\Pic(X_K)\cong \Pic(X_{\bar k}).
\end{equation}
The corresponding Leray spectral sequence gives rise to an
exact sequence
\begin{align}
\label{importantexaseq}
\begin{split}
0\to \ker(\Br(k)\to\Br(K))\to
\ker(\Br(X)\to \Br(X_K)) & \stackrel{\lambda}\to
H^1(G,\Pic(X_K)) \\
& \qquad\qquad \to H^3(G,K^*).
\end{split}
\end{align}
The inflation map
\begin{equation}
\label{inflation}
H^1(G,\Pic(X_K))\to H^1(\Gal(\bar k/k),\Pic(X_{\bar k})),
\end{equation}
is an isomorphism by the
Hochschild--Serre spectral sequence of group cohomology
(the inflation map is injective, and the cokernel maps into
$H^1$ of a torsion-free module with trivial action of a profinite group,
which is trivial).

Summarizing, we have

\begin{prop}
\label{prop:easy}
The composition of the inflation map \eqref{inflation} and the isomorphism
\eqref{BralgisH1} is an isomorphism
\begin{equation}
\label{eqn:easy}
\ker(\Br(X)\to\Br(X_{\bar k}))/\Br(k)\cong H^1(G,\Pic(X_K)).
\end{equation}
\end{prop}

In particular, the set $X(\A_k)^{\mathrm{Br.alg}}$
in \eqref{defXBralg}
is determined by finitely many
coset representatives $A\in \ker(\Br(X)\to \Br(X_{\bar k}))$, i.e.,
coset representatives of generators of the finite group
\eqref{eqn:easy}.
Our main theorem is

\begin{theo}
\label{iscomputable}
If $X$ is as in Assumption \ref{assumption1} then
$X(\A_k)^{\mathrm{Br.alg}}$ is effectively computable.
\end{theo}

\begin{rema}
\label{isallbr}
In this paper, we focus exclusively on
the algebraic part of the Brauer group and on effectively computing
$X(\A_k)^{\mathrm{Br.alg}}$.
Whenever $\Br(X_{\bar k})=0$, we have
$$\Br(X)/\Br(k)\cong H^1(G,\Pic(X_K))$$
by \eqref{BralgisH1} and \eqref{inflation}, and by the definition
\eqref{defXBralg} we have
$$X(\A_k)^{\mathrm{Br.alg}}=X(\A_k)^\Br.$$
The condition $\Br(X_{\bar k})=0$ holds automatically in either of
the following cases:
\begin{itemize}
\item[(i)] $X_{\bar k}$ is rational (by birational invariance
of $\Br(X)$ for smooth projective $X$ over a field of characteristic
zero \cite[(III.7.4)]{gb}), or
\item[(ii)] $X$ is Fano (i.e.,
has ample anticanonical divisor)
and $\dim X \le 3$ \cite{ip}.
\end{itemize}
An example of a computation of a Brauer--Manin obstruction
based on \emph{transcendental} (i.e., non-algebraic)
elements of the Brauer group is given in \cite{wittenberg}.
\end{rema}

\begin{exam}
Let $X$ be a del Pezzo surface of degree $3$ given by an equation
of diagonal form
\begin{equation}
\label{eq:cubic}
ax^3+by^3+cz^3+dt^3=0,
\end{equation}
where $a$, $b$, $c$, $d$ are nonzero integers.
Then $X(\Q)\ne\emptyset\Leftrightarrow X(k)\ne\emptyset$ where
$k=\Q(e^{2\pi i/3})$.
A splitting field is
$$K=k(\sqrt[3]{b/a},\sqrt[3]{c/a},\sqrt[3]{d/a}).$$
Specifically, all 27 exceptional curves of $X_{\bar k}$ are defined over $K$.
The Galois action of $G=(\Z/3\Z)^3$ on $\Pic(X_K)=\Z^7$ can be
computed.
By \cite{ctks}, the result is:
\begin{equation}
\label{eq:cubicH1}
H^1(G,\Pic(X_K))=\left\{
\begin{array}{cl}
0 & \text{if one of $ab/cd$, $ac/bd$, $ad/bc$ is a cube,}\\
(\Z/3\Z)^2 & \text{if exactly $3$ of $a/b$, $a/c$, $\ldots$, $c/d$ are cubes,}\\
\Z/3\Z & \text{otherwise.}
\end{array}
\right.
\end{equation}
\end{exam}

\section{Descent for divisors}
\label{sec:descentdivisors}
In this section we explain how to make effective the
isomorphism
\begin{equation}
\label{eq:galoisinvarpic}
\Pic(X)\cong \Pic(X_{\bar k})^{\Gal(\bar k/k)}
\end{equation}
from Remark \ref{remaiso}.
While this is not a logical requirement for the proof of
Theorem \ref{iscomputable}
(in fact it relies on
results from
Sections \ref{sec:2cocycle} and \ref{sec:eff}), it is important
in practical applications.
More precisely, given the data of Assumption \ref{assumption1},
let $H$ be a subgroup of $G$, with corresponding intermediate field
extension $K_0=K^H$.
Then \eqref{eq:galoisinvarpic} applied to $X_{K_0}$
yields divisors on $X_{K_0}$ representing elements of
$\Pic(X_{\bar k})^H$.
In this way, effective implementation of \eqref{eq:galoisinvarpic}
sometimes allows $K$ to be replaced by a smaller splitting field,
or at least a smaller extension over which a sufficiently interesting
submodule of $\Pic(X_{\bar k})$ is
defined.

We start with an example, and then explain how to carry this out in general.

\begin{exam}
Let $k=\Q(\zeta)$, where $\zeta=e^{2 \pi i/3}$,
and let $X$ be the diagonal cubic surface
\eqref{eq:cubic}.
Recall, we can take $K=k(\sqrt[3]{b/a},\sqrt[3]{c/a},\sqrt[3]{d/a})$ and
$G=(\Z/3)^3$.
Assume the coefficients $a$, $b$, $c$, and $d$ to be generic, so that
\begin{equation}
\label{eq:H1picX3}
H^1(G,\Pic(X_{\bar k}))=\Z/3\Z
\end{equation}
(see \eqref{eq:cubicH1}).
If we consider the subfield
$$K_0=k(\sqrt[3]{ad/bc})$$
then a computation reveals that the inflation map of Galois cohomology
\begin{equation}
\label{eq:infla}
H^1(\Gal(K_0/k),\Pic(X_{K_0}))\to H^1(G,\Pic(X_{\bar k}))
\end{equation}
is an isomorphism.
Concretely,
$\Pic(X_{K_0})\cong \Z\cdot(-\omega_X)\oplus M$ where $M$ is a rank $2$ module
with nontrivial action of $\Gal(K_0/K)\cong \Z/3\Z$.
The Galois-invariant combinations of exceptional lines on $X_{\bar k}$
generate only an index $3$ subgroup of $M$.
An additional generator of $M$ is the class of the following divisor
\begin{equation}
\label{eq:divisoroncubic}
D=D'-D'',
\quad
D':
\left\{
\begin{array}{l}
x+\zeta^2\sqrt[3]{b/a}\,y=0\\
z+\sqrt[3]{d/c}\,t=0
\end{array}
\right.
\quad
D'':
\left\{
\begin{array}{l}
x+\sqrt[3]{b/a}\,y=0\\
z+\zeta^2\sqrt[3]{d/c}\,t=0
\end{array}
\right.
\end{equation}
for which we do not, \emph{a priori}, have a representative defined over $K_0$.
Notice that $D$ is defined over
$$K_1=K_0(\sqrt[3]{b/a})=k(\sqrt[3]{b/a},\sqrt[3]{d/c}).$$
We define
$$\cL_1=\cO_{X_{K_1}}(D).$$

To make the isomorphism
\eqref{eq:galoisinvarpic} effective,
we apply
the following strategy:
We use \emph{descent} to produce a line bundle $\cL_0$
defined over $K_0$
having the desired class in the Picard group.
A rational section of $\cL_0$ defined over $K_0$ will produce
the required cycle.

The theory of descent is a machinery for patching,
i.e., the construction of a global object from local data
(see \cite[exp.\ VIII]{sga1}).
In this case, the local data consists of the line bundle
$\cL_1$ on $X_{K_1}$,
together with isomorphisms of
$\cL_1$ with its translates under $\Gal(K_1/K_0)$.
The isomorphisms which we must supply need to satisfy
a compatibility condition called the \emph{cocycle condition}.

For the sake of illustration,
we carry this out for the special choice of coefficients
$$a=5,\qquad b=9,\qquad c=10,\qquad d=12.$$
These coefficients are those of the famous example of Cassels and Guy
of a cubic surface which violates the Hasse principle \cite{cg}.
Then
$$K_0=k(\sqrt[3]{2/3}),\qquad\text{and}\qquad K_1=K_0(\sqrt[3]{9/5}).$$
Let $\rho$ denote a generator of $\Gal(K_1/K_0)\cong \Z/3\Z$.
To apply descent, we need to produce
an isomorphism
$\varphi\colon \cL_1 \to {}^\rho\!\cL_1$ satisfying the
cocycle condition
${}^{\rho\rho}\!\varphi\circ {}^\rho\!\varphi\circ\varphi=\mathrm{id}$;
then the descent machinery produces a line bundle $\cL_0$ on $X_{K_0}$.
This is easy to do because $X$ has a $K_0$-point, e.g.,
$$p:=(3,1,0,-\sqrt[3]{12}).$$
If we require
$\varphi$ to act as identity on the fiber of
$\cL_1$ over $p$,
then $\varphi$ is uniquely determined (since it is unique up to scale):
\begin{equation}
\label{eqn:phiinexa}
\varphi=-\frac{1}{2}(1+\sqrt[3]{15})
\frac{z+\sqrt[3]{6/5}\,t}{x+\sqrt[3]{9/5}\,y}.
\end{equation}
Starting with the function $1$ (viewed as a rational
section of $\cL_1$ or any of its Galois translates),
the rational section
$1+{}^{\rho\rho}\!\varphi(1)+
{}^{\rho\rho}\!\varphi\circ{}^{\rho}\!\varphi(1)$
is compatible with the Galois action, hence
descends to a rational section $s$ of $\cL_0$.
The divisor associated with the rational section $s$, which
must have the same class in the Picard group as $D$, is
$$C-(L+{}^\rho\!L+{}^{\rho\rho}\!L)$$
where $C$ is the cubic curve on $X_{K_0}$ defined by the equations
\begin{gather}
\label{eq:C}
\begin{split}
2x^2-6xy-xz+3\zeta\sqrt[3]{2/3}\,xt+3yz-9\zeta\sqrt[3]{2/3}\,yt+8z^2=0,
\,\,\,\,\,\,\,\,\,\,\,\\
4x^2-2xz-6\sqrt[3]{2/3}\,xt-6\zeta^2yz+z^2+3\sqrt[3]{2/3}\,zt
+9\sqrt[3]{4/9}\,t^2=0,\,\,\,\,\,\\
-2xy-5\zeta xz-\zeta^2\sqrt[3]{2/3}\,xt+6y^2-\zeta yz+3\zeta^2\sqrt[3]{2/3}\,yt
-8\sqrt[3]{2/3}\,zt=0,
\end{split}
\end{gather}
and where $L$ is the exceptional curve defined by
$x+\sqrt[3]{9/5} y=0$ and $z+\zeta^2\sqrt[3]{6/5} t=0$.
Now $\Pic(X_{K_0})$ is generated by an anticanonical divisor and
$C$ together with its translates under $\Gal(K_0/k)$.
\end{exam}

We return to the general setting.
Let $X$ be as in
Assumption \ref{assumption1}.
The machinery of descent associates, to a
vector bundle $\widetilde{\cE}$ on $X_K$
(or more generally a quasi-coherent sheaf of $\cO_{X_K}$-modules)
together with a collection of isomorphisms
$\varphi_g\colon \widetilde{\cE}\to a_g^*\widetilde{\cE}$
(for all $g\in G$) satisfying the cocycle condition
\begin{equation}
\label{eq:cocycle}
\varphi_{gh}=({}^g\!\varphi_h)\circ \varphi_g
\end{equation}
(for all $g$, $h\in G$)
a vector bundle (or quasi-coherent sheaf) $\cE$ on $X$ together with
an isomorphism $\xi\colon \cE_K\to \widetilde{\cE}$, which is
compatible with the $\varphi_g$.
Here ${}^g\!\varphi_h$ denotes $a_g^*\varphi_h$, where
$a_g$ is the automorphism of $X_K$ induced by $g\in G$.
(The reversed order of the composition on the right-hand side of
\eqref{eq:cocycle} is accounted for by our convention, in which the
Galois action induces a \emph{right} action of $G$ on $X$.)
The $\cE$ and $\xi$ that are produced by descent
are unique up to canonical isomorphism.

Let $D\subset X_K$ be a divisor (given by
equations)
whose class $[D]\in\Pic(X_{\bar k})$ is invariant under
$G=\Gal(K/k)$.
By Proposition \ref{DEprop},
if $D'\subset X_K$ is a divisor with
$[D']=[D]$ in $\Pic(X_{\bar k})$, then there is an effective procedure
to construct a rational function in $K(X)^*$ whose associated divisor is
$D-D'$.
Multiplication by this function is then an explicit isomorphism
$\cO_{X_K}(D)\to \cO_{X_K}(D')$.
For each $g\in G$, let
$$\varphi_g\colon \cO_{X_K}(D)\to \cO_{X_K}({}^g\!D)$$
be such an isomorphism.

Each isomorphism $\varphi_g$ is uniquely determined up to a
multiplicative constant.
We can characterize whether it is possible to modify
each isomorphism by a multiplicative constant, in order to satisfy
\eqref{eq:cocycle}.
The obstruction to \eqref{eq:cocycle} is the class in
$H^2(G,k^*)=\ker(\Br(k)\to\Br(K))$ of the following
$K^*$-valued $2$-cocycle
$(\gamma_{g,h})$:
$$\gamma_{g,h}:=\varphi_{gh}^{-1}\circ {}^g\!\varphi_h\circ \varphi_g.$$
As we have seen in the Example, the condition
$X(k)\ne\emptyset$ is sufficient for the obstruction to vanish.
In general we do not know whether $X(k)$ is empty.
But we have, by assumption,
$X(k_v)\ne\emptyset$ for all completions $k_v$ of $k$.
So the obstruction vanishes upon base change to any completion of $k$;
hence by the exact sequence \eqref{Brkseq} (the Hasse principle for the
Brauer group of a number field)
the obstruction indeed vanishes.
In other words, $(\gamma_{g,h})$ must be a $2$-coboundary.
There is an effective algorithm to express $(\gamma_{g,h})$ as the
coboundary of a $1$-cochain with values in $K^*$ (see the
proof of Proposition \ref{effectiveH3}).
So the $\varphi_g$ can be modified, using these multiplicative factors,
in order to satisfy \eqref{eq:cocycle}.
Applying descent, we obtain a line bundle $\cL$ on $X$,
such that $\cL$ is a representative of
$[D]\in\Pic(X_{\bar k})$.

It remains to express the line bundle $\cL$ (which is determined by
means of descent) explicitly as the class of a divisor on $X$.
Given any rational section of $\cL$, the associated divisor
(of zeros minus poles,
with respect to local trivializations of $\cL$)
will be defined over $k$
and will have class in the Picard group equal to $[D]$.
The theory of descent also includes descent for sections:
a rational section of $\cL$ is determined uniquely by a tuple
of rational sections of
$\cO_{X_K}({}^g\!D)$, for each $g$, that are compatible with the $\varphi_g$.
We obtain such a collection of sections from a single rational
section of
$\cO_{X_K}(D)$, by means of the
``trace'' operation: we translate the given section by all the elements of $G$
and form the sum of the translates.
So it suffices to exhibit a rational section of
$\cO_{X_K}(D)$ having nontrivial trace.
Let $x\in X_K$ be a closed point (not necessarily a $K$-rational point),
not lying in $D$ or in any of its Galois translates,
with Galois orbit $O(x)$ and image $y\in X$.
The extension from the residue field of $y$ to the coordinate ring
of $O(x)$, with $G$-action on the latter, can be calculated explicitly.
An element of the residue field of $x$
with nontrivial trace can be produced and
lifted to an element of
$H^0(X_K\smallsetminus\bigcup_{g\in G}{}^g\!D,\cO_{X_K}(D))$,
also with nontrivial trace.

\section{Computing the obstruction}
\label{sec:procedure}
In this section we explain the main steps of the computation
of $X(\A_k)^{\mathrm{Br.alg}}$
in terms of the data of Assumption \ref{assumption1}.
The details will be provided in subsequent sections, completing the
proof of Theorem \ref{iscomputable}.

We first need some additional notation. Put
\begin{equation}
\label{defU}
U = X \smallsetminus \bigcup_{i=1}^m D_i.
\end{equation}
Note that $U$ is defined over $k$.
We have an exact sequence
\begin{equation}
\label{Rseq}
0\to K^*\to \cO(U_K)^*\to R\to 0.
\end{equation}

\bigskip
\noindent
\emph{Step 1.} Compute the Galois cohomology group $H^1(G,\Pic(X_K))$,
and exhibit (finitely many) explicit $1$-cocycle representatives of
generators.

\bigskip
For \emph{each} generator (with cocycle representative) we implement the following.

\bigskip
\noindent
\emph{Step 2.} Apply the connecting homomorphism
of the cohomology exact sequence
of \eqref{picseq} to the cocycle representative of the generator from
Step 1 to obtain a $2$-cocycle representative of the corresponding element
in $H^2(G,R)$.

\bigskip
\noindent
\emph{Step 3.} Extending $K$ if necessary,
kill the obstruction in $H^3(G,K^*)$ to lifting
the element
of $H^2(G,R)$ to an element
$B\in H^2(G,\cO(U_K)^*)$;
carry out the lifting explicitly
on the cocycle level.

\bigskip
The element $B\in H^2(G,\cO(U_K)^*)$ will be
the restriction to $U$
of an element $A\in \Br(X)$.
More precisely the Leray spectral sequence
$$\mathsf{E}_2^{p,q}=H^p(G,H^q(U_K,\G_m))\Rightarrow H^{p+q}(U,\G_m)$$
induces a map
\begin{equation}
\label{eqn:LerayU}
H^2(G,\cO(U_K)^*)\to \Br(U).
\end{equation}
This map sends the class of a
$2$-cocycle with values in $\cO(U_K)^*$
to the element of $\Br(U)$ represented by the same cocycle, viewed now
as a $\check{\rm{C}}$ech cocycle for the covering $U_K\to U$.
Then Proposition \ref{restriction}, below, exhibits the required $A\in\Br(X)$.
The $\Br(k)$-coset $\widetilde{A}$ of $A$ will be one of finitely many generators of
$\ker(\Br(X)\to \Br(X_{\bar k}))/\Br(k)$.
Define
$$X(\A_k)^A=\{\,(x_v)\in X(\A_k)\,|\,\sum_v \inv_v(A(x_v))=0\,\}.$$
Then
$$X(\A_k)^{\mathrm{Br.alg}}=\bigcap_{A} X(\A_k)^A,$$
where $A$ runs over the finite set of representatives.

\bigskip
\noindent
\emph{Step 4.} Obtain from the set $\mathcal{D}^1:=\{D_i\}$,
new finite collections $\mathcal{D}^2$, $\ldots$, $\mathcal{D}^r$
of geometric cycles, such that the corresponding complements
$U=:U^1$, $U^2$, $\ldots$, $U^r$ form an open covering of $X$.
Repeat Steps 2 and 3 for each $\mathcal{D}^j$ to obtain
$B^j\in H^2(G,\cO(U^j_K)^*)$ such that the restriction of $A$
to $U^j$ is equal to $B^j$, modulo $\Br(k)$.

\bigskip
\noindent
\emph{Step 5.} (Calibration) Compute $I^j\in \Br(k)$ such that
$$B^j+I^j=A|_{U^j}$$
in $\Br(U^j)$.

\bigskip
\noindent
\emph{Step 6.} Compute the local invariants $\inv_v(A(x_v))$ for all
$v$ and all $x_v\in X(k_v)$.

\section{$2$-cocycle representatives}
\label{sec:2cocycle}
In this section we carry out Steps 1 through 3 outlined above.
We obtain $2$-cocycle representatives
for the classes $B\in H^2(G,\cO(U_K)^*)$ from Section \ref{sec:procedure}.

Step 1, the computation of $H^1(G,\Pic(X_K))$, is implemented
in standard computer algebra packages, e.g., \texttt{Magma}.
Indeed, by Assumption \ref{assumption1}, $\Pic(X_K)$ with its
Galois action is known.
The output is a presentation of $H^1(G,\Pic(X_K))$ as a finite
abelian group, with $1$-cocycle representatives of a set of generators.

Steps 2 and 3 produce lifts of a generator to $\Br(X)$, via the map
$\lambda$ of the sequence \eqref{importantexaseq}.
All computations are done on the level of cocycle representatives.
The obstructions in $H^3(G,K^*)$ to producing the lift are killed
by enlarging $K$, if necessary.

Let $\widetilde{B}\in H^1(G,\Pic(X_K))$ be one of the generators, with
$1$-cocycle representative
$\tilde\beta$.
Concretely, $\tilde\beta$ is a tuple of elements of $\Pic(X_K)$,
satisfying a cocycle condition.

Combining the cohomology exact sequences coming from
\eqref{picseq} and \eqref{Rseq} and a portion of the
exact sequence \eqref{importantexaseq} we obtain a diagram
\begin{equation}
\label{diagram1}
\begin{split}
\xymatrix@C=12pt@R=15pt{
&&&{\Br(U)}\\
&{\ker(\Br(X)\to\Br(X_K))}\ar[urr]\ar[d]^\lambda&
H^2(G,\cO(U_K)^*)\ar[ur]\ar[d]^\mu & \\
0 \ar[r] &
H^1(G,\Pic(X_K)) \ar[r]^(0.57)\delta \ar[dr]_{\varepsilon} &
H^2(G,R) \ar[r] \ar[d]^\nu &
H^2(G,\bigoplus_{i=1}^m \Z\cdot[D_i]) \\
&&H^3(G,K^*)&
}
\end{split}
\end{equation}
In this diagram, we have used the fact that $\bigoplus_{i=1}^m \Z\cdot[D_i]$
is a permutation module, and hence its first cohomology vanishes.
The maps to $\Br(U)$ are the restriction map from $\Br(X)$ and
the map from $H^2(G,\cO(U_K)^*)$ of \eqref{eqn:LerayU}, respectively.

\begin{prop}
\label{restriction}
Let $X$ be as in Assumption \ref{assumption1},
let $\widetilde{B}\in H^1(G,\Pic(X_K))$.
\begin{itemize}
\item[(i)] The map $\varepsilon$ in \eqref{diagram1} is
the rightmost map in \eqref{importantexaseq}.
\item[(ii)] Suppose that
$\varepsilon(\widetilde{B})=0$.
Let $A\in \ker(\Br(X)\to\Br(X_K))$ be a lift of $\widetilde{B}$ by
the map $\lambda$ and
$B\in H^2(G,\cO(U_K)^*)$ a lift of
$\delta(\widetilde{B})$ by the map $\mu$ in \eqref{diagram1}.
Then the images of $A$ and $B$ in $\Br(U)$ via the maps in
\eqref{diagram1} are equal modulo $\Br(k)$.
\end{itemize}
\end{prop}

\begin{prop}
\label{splitting}
With notation as above, there is an effective construction
of a splitting $\sigma\colon R\to \cO(U_K)^*$ of
the sequence \eqref{Rseq}.
\end{prop}

\begin{prop}
\label{effectiveH3}
Let $k$ be a number field, and $K/k$ a finite Galois extension
with Galois group $G=\Gal(K/k)$.
Let $\kappa$ be a cocycle representative of an
element in $H^3(G,K^*)$.
There is an effective algorithm to determine whether $\kappa$ is trivial
in $H^3(G,K^*)$.
If so, we can effectively produce a lift of $\kappa$ to a
$2$-cochain via the coboundary map.
Otherwise,
we can effectively produce
a finite extension $L/K$, with $L$ Galois over $k$,
such that the inflation of $\kappa$ to $H^3(\Gal(L/k),L^*)$ is trivial.
\end{prop}

It is straightforward to compute $\delta(\tilde\beta)$.
Now $\varepsilon(\tilde\beta)=\nu(\delta(\tilde\beta))$,
and $\nu(\delta(\tilde\beta))$ is computed using
the splitting $\sigma\colon R\to \cO(U_K)^*$ of
Proposition \ref{splitting}.
Proposition \ref{effectiveH3} supplies an extension $L$ of $K$
killing the class of $\varepsilon(\tilde\beta)$ in $H^3(G,K^*)$.
Replacing $K$ by $L$, we now invoke Proposition \ref{effectiveH3}
to produce a $2$-cochain $\eta$ with values in $K^*$, such that
$$\varepsilon(\tilde\beta)=\partial(\eta)$$
Now put
$$\beta:=\frac{\sigma(\delta(\tilde\beta))}{\eta}.$$
Then $\beta$ is a $2$-cocycle with values in $\cO(U_K)^*$,
such that the class of
$\beta$ is a lift via $\mu$ of the
class $\delta(\widetilde{B})$.
We let $B\in H^2(G,\cO(U_K)^*)$ denote the class of $\beta$.

By the exact sequence
\eqref{importantexaseq} and Proposition \ref{restriction} (i)
the class $\widetilde{B}\in H^1(G,\Pic(X_K))$ lifts via $\lambda$ to a class
\begin{equation}
\label{Ai}
A\in \ker(\Br(X)\to\Br(X_K)),
\end{equation}
defined up to an element of $\Br(k)$.
By Proposition \ref{restriction} (ii),
$A$ can be chosen so that
the image of $B$ in $\Br(U)$ is equal to the restriction $A|_U$.
Thus, in \eqref{Ai}, we have a Brauer group element $A$ whose restriction
to $U$ is known explicitly, by means of the $2$-cocycle $\beta$.

\begin{proof}[Proof of Proposition \ref{restriction}]
We follow \cite{cts2}.
Let $j\colon U\to X$ be the inclusion, with complement
$D=\bigcup_{j=1}^m D_i$.
There is an exact sequence of \'etale sheaves on $X$
\begin{equation}
\label{GmKZ}
0\to \G_m\to j_*(\G_m|_U)\to \mathcal{Z}^1_D\to 0
\end{equation}
where on the right is the sheaf of divisors on $X$ with support on $D$.
Evaluating global sections on $X_K$,
we obtain the exact sequence of $G$-modules
\begin{equation}
\label{2termextn}
0\to K^*\to \cO(U_K)^*\to \bigoplus_{i=1}^m \Z\cdot[D_i]\to
\Pic(X_K)\to 0.
\end{equation}

The sequence \eqref{GmKZ} provides a resolution of the sheaf $\G_m$ on $X$.
Taking $(\mathcal{I}_\bullet,\mathsf{d})$ to be a resolution of $\G_m$ by injective
sheaves on $X$, we know that there exists a morphism of resolutions
from the resolution \eqref{GmKZ} to $\mathcal{I}_\bullet$.
Applying the equivariant global section functor, we get a morphism
of four-term exact sequences
\begin{equation}
\label{another2termextn}
\begin{split}
\xymatrix@R=20pt@C=16pt{
K^*\ar@{^{(}->}[r]\ar@{=}[d]
&{\cO(U_K)^*}\ar[r]^{\dv}
\ar[d]^{\psi}&{\bigoplus_{i=1}^m \Z\cdot[D_i]}\ar@{->>}[r]\ar[d]^{\varphi}
&{\Pic(X_K)}\ar@{=}[d]\\
K^*\ar@{^{(}->}[r]&H^0(X_K,\mathcal{I}_0)\ar[r]^(0.30){\mathsf{d}}&
{\ker(H^0(X_K,\mathcal{I}_1)\to H^0(X_K,\mathcal{I}_2))}\ar@{->>}[r]
&{\Pic(X_K)}
}
\end{split}
\end{equation}
where the first and last maps are identity maps.

The four-term exact sequence \eqref{2termextn}
is the amalgamation of two short exact sequences,
and the map $\varepsilon$ in the diagram \eqref{diagram1} is the
composition of the connecting homomorphisms of the two long exact sequences in
cohomology.
Spectral sequence machinery shows that the edge map
$H^1(G,\Pic(X_K))\to H^3(G,K^*)$ of the sequence \eqref{importantexaseq}
is equal to a similar composition of connecting homomorphisms coming from the
four-term exact sequence of the bottom line of \eqref{another2termextn}.
Because there is a map between these sequences inducing identity maps on
the first and last terms, the morphism $\varepsilon$ is equal to the
edge map of \eqref{importantexaseq}.
This establishes part (i).

For (ii),
let $\widetilde{B}\in H^1(G,\Pic(X_K))$ be given, represented by the $1$-cocycle
$\tilde\beta$.
Computing $\delta(\widetilde{B})$ involves lifting $\tilde\beta$ to
a $1$-cochain $\tilde\gamma$ with values in $\bigoplus_{i=1}^m \Z\cdot[D_i]$.
Now $\partial(\tilde\gamma)$ is a tuple of divisors
rationally equivalent to zero, hence
$\partial(\tilde\gamma)=\dv(\beta)$
for some $2$-cochain $\beta$ with values in $\cO(U_K)^*$.
Moreover $\partial(\beta)$ is a $3$-cocycle representative of
$\varepsilon(\widetilde{B})$ which by hypothesis vanishes.
So, by adjusting $\beta$ by a $K^*$-valued $2$-cochain
we may arrange that $\partial(\beta)=0$.
Then $\beta$ is a $2$-cocycle representative of a class
$B\in H^2(G,\cO(U_K)^*)$, such that $\mu(B)=\delta(\widetilde{B})$.

We may identify $\Br(X)$ with the second cohomology of the total complex
of the term $\mathsf{E}_0^{p,q}=\mathsf{C}^p(H^0(X_K,\mathcal{I}_q))$ of the Leray spectral
sequence.
Now $A\in \Br(X)$ is represented by a cocycle of the total complex
$$(\alpha_0,\alpha_1,\alpha_2)\in \mathsf{C}^0(H^0(X_K,\mathcal{I}_2))\times
\mathsf{C}^1(H^0(X_K,\mathcal{I}_1))\times \mathsf{C}^2(H^0(X_K,\mathcal{I}_0)).$$
Since $A\in \ker(\Br(X)\to\Br(X_K))$, we may suppose that
$\alpha_0=0$.
The condition to be a cocycle is now
$$\mathsf{d}(\alpha_1)=0,\qquad\qquad
\partial(\alpha_1)=\mathsf{d}(\alpha_2),\qquad\qquad \partial(\alpha_2)=0.$$
The cocycle representative $(\alpha_0,\alpha_1,\alpha_2)$
may be replaced by an equivalent representative with
$\alpha_1=\varphi(\tilde\gamma)$.
Then
$\mathsf{d}(\psi(\beta))=\varphi(\dv(\beta))=
\partial(\varphi(\tilde\gamma))=\mathsf{d}(\alpha_2)$.
Hence the image of $A-B$ in $\Br(U)$ is represented by a
$2$-cocycle with values in $\ker(\mathsf{d})=K^*$.
\end{proof}

\begin{proof}[Proof of Proposition \ref{splitting}]
The map to $R$ in \eqref{Rseq}, for which we wish to find a splitting,
is the divisor map from rational functions on $X$,
regular and nonvanishing on $U$, to divisors with support outside $U$.
So it suffices to solve the problem,
given effective divisors $D$ and $E$ on a smooth
projective variety $X$ over $k$
(all given explicitly by equations), with $[D]=[E]$ in $\Pic(X)$, to
produce effectively a rational function in $k(X)^*$ whose divisor is
$D-E$.
This is the content of Proposition \ref{DEprop},
given in the next section.
\end{proof}

\begin{proof}[Proof of Proposition \ref{effectiveH3}]
By \cite{bsd} Theorem 3, there is an effective method to test whether
$\kappa=0$ in $H^3(G,K^*)$, and to produce a lift to a $2$-cochain if this
is the case.
The method is to produce, effectively, a finite set of primes $S$
such that $\kappa=0$ if and only if the given cocycle is the coboundary
of a $2$-cochain taking values in the $S$-units of $K$.
(The same argument applies to test for triviality of an $i$-cocycle,
for any $i$, and to produce an $(i-1)$-cochain in case it is trivial.)

If $\kappa\ne 0$ in $H^3(G,K^*)$, then there exist
cyclic extensions $\ell$ of $k$ such that
$L:=\ell K$ satisfies $\kappa\in\ker(H^3(G,K^*)\to H^3(\Gal(L/k),L^*))$.
For instance, let $q$ be a prime not dividing the discriminant
$\disc(K/\Q)$ such that $q\equiv 1\pmod{n}$.
Then $[L:k]=(q-1)n$ where $n=[K:k]$.
By the Chebotarev density theorem, there exists
some prime ideal $\mathfrak{p}$ in $k$ (which we do not need explicitly)
which remains inert in the cyclic extension $\ell$ of $k$.
Then the local degree $n_\mathfrak{p}$ of $\mathfrak{p}$ in $L$ must
be a multiple of $q-1$.
Therefore the inflation map
$H^3(G,K^*)\to H^3(\Gal(L/k),L^*)$ is trivial
(cf.\ \cite{at}, Section 7.4).
\end{proof}

\begin{exam}
For the cubic surface over $k=\Q(\zeta)$, $\zeta=e^{2\pi i/3}$, defined by
$$5x^3+9y^3+10z^3+12t^3=0,$$
we have found that already for
the cyclic degree $3$ extension $K_0=k(\sqrt[3]{2/3})$ we have
$$H^1(\Gal(K_0/k),\Pic(X_{K_0}))\cong H^1(G,\Pic(X_{\bar k}))=\Z/3\Z.$$
Let $\tau$ denote the generator of $\Gal(K_0/k)\cong \Z/3\Z$ which
sends $\sqrt[3]{2/3}$ to $\zeta\sqrt[3]{2/3}$.
Cohomology $H^i(\Gal(K_0/k),M)$ can be computed by means of the resolution
\begin{equation}
\label{eq:cyclicgroupcoho}
\mathsf{C}^\bullet(M):\qquad 0\longrightarrow
M\stackrel{\Delta_\tau}\longrightarrow
M\stackrel{N_\tau}\longrightarrow
M\stackrel{\Delta_\tau}\longrightarrow\cdots
\end{equation}
where the maps alternate between $\Delta_\tau:=\mathrm{id}-\tau$
and $N_\tau:=\mathrm{id}+\tau+\tau^2$.
The group $H^1(\Gal(K_0/k),\Pic(X_{K_0}))$ is thus identified with
$\ker(N_\tau)/\im(\Delta_\tau)$, and this group
(which is cyclic of order $3$) is generated by the class of
\begin{equation}
\label{eqn:theclass}
[C]+\omega_X
\end{equation}
where $C$ is the curve in $X_{K_0}$ given in \eqref{eq:C}.
This is Step 1.

Choose an anticanonical divisor $H\subset X$ defined, say, by $x=0$,
so that $C-H$ is a divisor in the class \eqref{eqn:theclass}.
Now $\{H,C,{}^\tau\!C,{}^{\tau\tau}\!C\}$ is a Galois-invariant
set of divisors generating $\Pic(X_{K_0})$.
For this set of divisors the sequence \eqref{picseq} becomes
\begin{equation}
\label{picseqexa}
0\to \Z\to \Z\cdot H\oplus \Z\cdot C\oplus \Z\cdot {}^\tau\!C
\oplus \Z\cdot {}^{\tau\tau}\!C \to \Pic(X_{K_0})\to 0.
\end{equation}

The connecting homomorphism induced by the resolution
\eqref{eq:cyclicgroupcoho} on the exact sequence \eqref{picseqexa}
sends the $1$-cocycle \eqref{eqn:theclass} to the $2$-cocycle
\begin{equation}
\label{eq:CtCttC}
C+{}^\tau\!C+{}^{\tau\tau}\!C-3H
\end{equation}
in $R$.
Notice that in this example $R$ is isomorphic to $\Z$ (the left-hand term in
\eqref{picseqexa}), and in fact \eqref{eq:CtCttC} is a generator.
We have completed Step 2.

We have $U=X\smallsetminus(H\cup C\cup {}^\tau\!C\cup {}^{\tau\tau}\!C)$.
A splitting $\sigma\colon R\to \cO(U_{K_0})^*$, required for Step 3,
sends the generator \eqref{eq:CtCttC} to
a rational function which vanishes on
$C\cup {}^\tau\!C\cup {}^{\tau\tau}\!C$
and has a pole of order $3$ along $H$.
The obstruction group $H^3(\Gal(K_0/k),K_0^*)$ is identified by
\eqref{eq:cyclicgroupcoho} with $H^1(\Gal(K_0/k),K_0^*)$,
which vanishes by Hilbert's Theorem 90.
So it is possible to lift the $2$-cocycle \eqref{eq:CtCttC}
to an element of $\cO(U_{K_0})^*$ invariant under $\tau$.
To carry this out, we use the
explicit equations \eqref{eq:C} for $C$
to produce directly a function in
$\cO(U)^*$ of the form
$f/x^3$, where $f$ is a polynomial (with coefficients in $k$)
whose divisor is $C+{}^\tau\!C+{}^{\tau\tau}\!C-3H$.
We obtain
\begin{align}
\label{eq:vanishonCtCttC}
\begin{split}
f=(2\zeta-2)x^3-3\zeta x^2y-8\zeta x^2z-{}&9\zeta^2xy^2+24\zeta xyz
+4\zeta xz^2+(-6\zeta-21)y^3\\
&{}-12\zeta yz^2+(-18\zeta-14)z^3+(4\zeta-4)t^3.
\end{split}
\end{align}
So we have completed Step 3.
The $2$-cocycle given by $f/x^3$ corresponds to a cyclic Azumaya
algebra for the extension $K_0$ of $k$ and the rational function
$f/x^3\in \cO(U)^*$.
The class $B\in \Br(U)$ of this Azumaya algebra
is the restriction of some $A\in \Br(X)$ such that
$A$ generates $\Br(X)/\Br(k)$.
\end{exam}

\section{Effectivity}
\label{sec:eff}

In this section, we present effectivity results concerning
ample line bundles and homogeneous ideals.
Then we show how to carry out Step 4.

Many results in effective algebraic geometry are based on
Gr\"obner bases.
There are effective algorithms to compute a Gr\"obner basis
of a homogeneous ideal $\mathcal{I}\subset k[x_0,\ldots,x_N]$,
which is
given by means of generators.
Based on this, there are effective algorithms (implemented in
computer algebra packages) to:
\begin{itemize}
\item test whether a given polynomial is in $\mathcal{I}$,
\item compute the saturation of $\mathcal{I}$,
\item compute the
primary decomposition of an ideal $\mathcal{I}$.
\end{itemize}
(See, e.g., \cite{clo}, \cite{kr}.)

\begin{lemm}
\label{lemmeffd}
Let $X\subset \PP^N$ be a projective variety, given by means of equations.
Denote by $\cO_X(1)$ the restriction of $\cO_{\PP^N}(1)$ to $X$.
Let $\cL$ be an arbitrary line bundle on $X$, presented
(as a coherent sheaf on $X$)
by means of homogeneous generators and relations.
Then there is, effectively, a positive integer $d_0$ such that for all
$d\ge d_0$ the line bundle
$\cL\otimes \cO_X(d)$ is ample and generated by global sections.
\end{lemm}

\begin{proof}
Choose finitely many sections $s_i\in H^0(X,\cO_X(1))$ such that
$X$ is covered by the open subsets $X_{s_i}$ where each section is
nonvanishing. From the presentation of $\cL$ we have generators
$t_{i,j}\in H^0(X_{s_i},\cL)$.
These give rise to a covering of the affine scheme $X_{s_i}$ by
open affines $X_{t_{i,j}}$.

There is, effectively, a positive integer $d_0$ such that
for all $d\ge d_0$, one has that
$t_{i,j}\otimes s_i^d$ is the restriction to $X_{s_i}$ of a
section $u_{i,j}\in H^0(X,\cL\otimes \cO_X(d))$
for each $i$ and $j$.
This is achieved following the proof of \cite[9.3.1(ii)]{egaI}.
The open subsets $X_{u_{i,j}}=X_{t_{i,j}}$ cover $X$, hence
$\cL\otimes \cO_X(d)$ is generated by global sections,
and by \cite[4.5.2($\mathrm{a}'$)]{egaII}, $\cL\otimes \cO_X(d)$ is ample.
\end{proof}

\begin{rema}
\label{remamatsusaka}
By an effective Matsusaka theorem \cite{siu}
there is, for given $d\ge d_0$, an effective bound $r_0$
such that
$(\cL\otimes \cO_X(d))^{\otimes r}$ is very ample for all $r\ge r_0$.
\end{rema}

%
%

\begin{prop}
\label{DEprop}
Given effective divisors $D$ and $E$ (by means of equations) on a smooth
projective variety $X$ (also defined by given equations) over $k$,
such that $[D]=[E]$ in $\Pic(X)$, there is an effective algorithm to
produce a rational function $h\in k(X)^*$ such that
$\dv(h)=D-E$.
\end{prop}

\begin{proof}
Let $H$ denote a hyperplane section of $X\subset \PP^N$.
By effective Serre vanishing \cite[Prop.\ 1]{bel}
we have
\begin{equation}
\label{H0surj}
H^0(\PP^N,\cO_{\PP^N}(d))\to H^0(X,dH)
\end{equation}
surjective for $d$ greater than some effective lower bound.
Combining this with the lower bound from Lemma \ref{lemmeffd} applied to
$\cL=\cO_X(-D)$,
we obtain $d$ such that \eqref{H0surj} is surjective and such that
we can produce
a homogeneous polynomial $f$ of degree $d$, not identically zero on $X$,
and vanishing on $D$.
So, $f$ vanishes on a cycle $D+D'$ for some cycle
$D'$ that can be determined effectively.

Now $E+D'$ is linearly equivalent to $dH$.
By surjectivity of \eqref{H0surj}, there exists
a homogeneous polynomial $g$ of degree $d$, not identically zero on $X$,
contained in the ideal of the sum of divisors $E+D'$.
The ideals can be computed effectively, so we can find $g$, and then
defining $h$ to be $f/g$ we have
$\dv(h)=D-E$.
\end{proof}

\begin{rema}
\label{remaavoidpoints}
Since the linear system $|dH-D|$ that arises in the proof of
Proposition \ref{splitting} is base point free, the divisor $D'$
in the proof can be chosen to avoid any given finite set of points
of $X$.
\end{rema}

\begin{rema}
\label{shrinkU}
The elements $f/g$ constructed in the proof of
Proposition \ref{splitting} do indeed lie in $\cO(U_K)^*$.
However, $f$ and $g$ individually vanish on the residual divisor $D'$
that appears in the proof.
For the purposes of effective computations, we would really require
alternative representations of this rational function, such that at
any point of $U$ there is some representative that is amenable to evaluation.
A way of avoiding this extra complication is to carry out the
following steps:
\begin{itemize}
\item[(i)] Fix a choice of hyperplane section $H$ of $X$ (in its given
projective embedding), which we assume to be defined by a linear polynomial
$h$ having coefficients in $k$ and included among
the cycles $D_i$.
\item[(ii)] Enlarge the collection of cycles $D_i$ by adding, for
each generator of $R$, the corresponding residual cycle $D'$ from the
proof of Proposition \ref{splitting}
along with its Galois translates.
\item[(iii)] For each additional cycle $D'$ (or Galois translate thereof),
add a new generator $D+D'-dH$ (or Galois translate thereof) to $R$.
\item[(iv)] For each additional generator to $R$, extend the definition
of the splitting $\sigma$ by sending $D+D'-dH$ to $f/h^d$ and each
Galois translate of $D+D'-dH$ to the corresponding Galois translate of $f/h^d$.
\end{itemize}
Having done this, $U$ will be the complement of the zero locus of
some homogeneous polynomial, and all the rational functions in the
image of $\sigma$ will be of the form $f/g$ where $f$ and $g$
are homogeneous polynomials nonvanishing on $U$.
\end{rema}

We now describe how to carry out Step 4.
For each divisor $D_i$, Lemma \ref{lemmeffd} applied to
$\cO_X(-D_i)$ yields $d$ such that $|dH-D_i|$ is base point free.
We may therefore choose a member of this linear system which avoids
any finite collection of points.
So we can obtain divisors $D'_i$, such that $D_i$ is in the
$\Z$-linear span of $D'_i$ and a hyperplane class, and such that
no irreducible component of any $D_i$ is contained in any of the $D'_j$.
We then take $\mathcal{D}^2$ to be the collection of $D'_i$ together
with a hyperplane class, the latter also chosen not to contain any
irreducible component of any of the $D_i$.

Now we repeat Steps 2 and 3, and also carry out steps (i)--(iv)
of Remark \ref{shrinkU} (which adds extra divisors to $\mathcal{D}^2$),
so that every rational function
that has been constructed from the divisors $D'_i$
has the property that its numerator and
denominator are both nonvanishing on $U^2$.
In carrying this out, we use Remark \ref{remaavoidpoints} to ensure that
$U^2$ contains the generic points of all the irreducible components of
the divisors in $\mathcal{D}^1$.

For $m=2$, $\ldots$, $\dim X$ we recursively
construct $\mathcal{D}^{m+1}$ and repeat Steps 2 and 3 and steps (i)--(iv)
of Remark \ref{shrinkU} for the divisors in $\mathcal{D}^{m+1}$,
making use of Remark \ref{remaavoidpoints} to ensure that $U^{m+1}$
contains all the generic points of $\bigcap_{i=1}^m (X\smallsetminus U^i)$.
It follows inductively that
$$\dim \Bigl(\bigcap_{i=1}^m (X\smallsetminus U^i)\Bigr)=\dim X-m.$$
At the end of the construction $X$ is covered by the open sets
$U^1$, $\ldots$, $U^{\dim X+1}$.

\section{Calibration}
\label{sec:calibration}
It is necessary to produce several representatives of a given
$\widetilde{A}\in \Br(X)/\Br(k)$ (Step 4)
and to calibrate these classes, i.e., to compute
the difference in $\Br(k)$ between two given representatives
in $\Br(X)$ of the $\Br(k)$-coset $\widetilde{A}$
(Step 5).

Let $\mathcal{D}=\{D_i\}$ and $\mathcal{D}'=\{D'_i\}$ be two
sets of divisors as at the end of Section \ref{sec:eff}.
Let $U=X\smallsetminus \bigcup D_i$ and
$U'=X\smallsetminus \bigcup D'_i$.
Let $\widetilde{B}$ be an element of $H^1(G,\Pic(X_K))$,
with cocycle representative $\tilde\beta$ (constructed in Step 1).
Consider the
$2$-cocycles $\beta$ and $\beta'$ resulting from Steps 2 and 3 applied to
$\mathcal{D}$ and $\mathcal{D}'$, respectively.
These are $2$-cocycles
with values in $\cO(U_K)^*$ and $\cO(U'_K)^*$, respectively.

The exact sequence \eqref{picseq} can be enlarged to an exact sequence
\begin{equation}
\label{enlargedpicseq}
0\to S\to \bigoplus_i \Z\cdot [D_i]\oplus \bigoplus_i \Z\cdot [D'_i]\to
\Pic(X_{\bar k})\to 0
\end{equation}
where $S$ contains $R$ as a direct summand.
There is also a sequence analogous to \eqref{Rseq},
\begin{equation}
\label{seqanalogous}
0\to K^*\to \cO(U_K\cap U'_K)^*\to S\to 0,
\end{equation}
and the splitting $\sigma$ of Proposition \ref{splitting} can be
extended to a splitting
$$\pi\colon S\to \cO(U_K\cap U'_K)^*.$$

Looking at Step 2, carried out using the collection of divisors $\mathcal{D}$ and
again using the collection of divisors $\mathcal{D}'$, we have the same
$1$-cocycle $\tilde\beta$ with values in $\Pic(X_{\bar k})$ mapped to the same
element of $H^2(G,S)$ represented by two different $2$-cocycles
with values in $S$.
So, these $2$-cocycles differ by a $2$-coboundary.
We apply $\pi$ (note that $\pi$ is compatible with the splitting
used in Step 3 for the collection of divisors $\mathcal{D}'$
only up to multiplicative constants),
to obtain
$$\beta = \iota'\,\beta'\,\partial(\pi(\theta))$$
for some $1$-cochain $\theta$ with values in $S$ and some
$2$-cochain $\iota'$ with values in $K^*$.
This $\iota'$ is a $2$-cocycle, and its associated class $I'\in \Br(k)$
satisfies
$$B'+I'=A.$$

\begin{exam}
In the case of the diagonal cubic surface
$$5x^3+9y^3+10z^3+12t^3=0,$$
we found that with $K_0=k(\sqrt[3]{2/3})$,
we have $\Pic(X_{K_0})$ generated by
the class of the hyperplane $H$ (given by $x=0$) and
the Galois orbit of the class of the cubic curve $C$ which we determined
in \eqref{eq:C} starting from a choice of $p\in X(K_0)$.
Additional choices of $K_0$-points $p'$, $\ldots$, lead to
further divisors having the same class in the Picard group as $C$.
For instance, the point
$$p':=(3\zeta,1,0,-\sqrt[3]{12})$$
gives rise to the curve $C'$ defined by
\begin{gather*}
2x^2-6\zeta xy-xz+3\zeta^2\sqrt[3]{2/3}\,xt+3\zeta yz-9\sqrt[3]{2/3}\,yt+8z^2=0,
\\
4x^2-2xz-6\zeta\sqrt[3]{2/3}\,xt-6yz+z^2+3\zeta\sqrt[3]{2/3}\,zt
+9\zeta^2\sqrt[3]{4/9}\,t^2=0,\\
-2xy-5xz-\zeta^2\sqrt[3]{2/3}\,xt+6\zeta y^2-\zeta yz+3\sqrt[3]{2/3}\,yt
-8\sqrt[3]{2/3}\,zt=0.
\end{gather*}
Applying Steps 2 and 3 to $C'$ we obtain
\begin{align}
\label{fprime}
\begin{split}
f'=(2\zeta+4)x^3-3\zeta x^2y-8x^2z-{}&9xy^2+24\zeta xyz+4xz^2+(21\zeta+15)y^3 \\
&{}-12\zeta yz^2+(14\zeta-4)z^3+(4\zeta+8)t^3\\
\end{split}
\end{align}
such that $\dv(f'/x^3)=C'+{}^\tau\!C'+{}^{\tau\tau}\!C'-3H$.

The sequence \eqref{enlargedpicseq} is
$$0\to \Z\oplus \Z^3\to \Z \oplus \Z^3 \oplus \Z^3\to \Pic(X_{K_0})\to 0$$
where in the middle the generators are $H$ and the Galois translates
of $C$ and $C'$.
The first generator of the rank $4$ module $S$ on the left is
the previously identified cycle $C+{}^\tau\!C+{}^{\tau\tau}\!C-3H$,
and the additional generators are
$$C'-C,\qquad
{}^\tau\!C'-{}^\tau\!C,\qquad
{}^{\tau\tau}\!C'-{}^{\tau\tau}\!C.
$$
We need to construct a splitting
$\pi\colon \Z\oplus \Z^3\to \cO(U_K\cap U'_K)^*$.
We already have the image $f/x^3$ for the first generator, with $f$ as in
\eqref{eq:vanishonCtCttC}.
Images of the other generators are readily constructed.
We send $C'-C$ to $g/f$, where
\begin{align*}
g&=(8\zeta + 16)x^3 + (-4\zeta - 8)x^2y - 2x^2z
+ (-2\zeta + 2)\sqrt[3]{2/3}\,x^2t
+ (2\zeta + 4)xyz\\
&{}+ (6\zeta - 6)\sqrt[3]{2/3}\,xyt + (-5\zeta + 6)xz^2 +
(\zeta - 1)\sqrt[3]{2/3}\,xzt + (-6\zeta - 3)\sqrt[3]{4/9}\,xt^2\\
&{}+ (12\zeta + 24)y^3 +
(-6\zeta + 6)y^2z + (-\zeta - 2)yz^2 + (-3\zeta + 3)\sqrt[3]{2/3}\,yzt\\
&{} + (18\zeta + 9)\sqrt[3]{4/9}\,yt^2
+ (16\zeta + 24)z^3 + (-8\zeta + 8)\sqrt[3]{2/3}\,z^2t + (16\zeta + 32)t^3
\end{align*}
and the remaining generators to the Galois conjugates of $g/f$.

The $2$-cocycles $C+{}^\tau\!C+{}^{\tau\tau}\!C-3H$ used to construct
$f$ and $C'+{}^\tau\!C'+{}^{\tau\tau}\!C'-3H$ used to construct $f'$
differ by a $2$-coboundary,
$$C+{}^\tau\!C+{}^{\tau\tau}\!C-C'-{}^\tau\!C'-{}^{\tau\tau}\!C'=
N_\tau(C-C').$$
We have
$$\pi(C-C')=\frac{f}{g}.$$
Consequently
$$\frac{f}{x^3}=\vartheta\frac{f'}{x^3}N_\tau\Bigl(\frac{f}{g}\Bigr)$$
for some constant $\vartheta$.
By an explicit computation, we find
$$\vartheta=\frac{\zeta}{4}.$$

Repeating everything with the point
$$p'':=(3,\zeta,0,-\sqrt[3]{12})$$
yields a curve $C''$ with explicit equations and a function
\begin{align}
\label{fprimeprime}
\begin{split}
f''=(-4\zeta-2)x^3-&{}3\zeta x^2y-8x^2z-9\zeta xy^2
  +(-24\zeta-24)xyz+4\zeta xz^2\\
&{}+(-15\zeta+6)y^3-12yz^2+(4\zeta+18)z^3
  +(-8\zeta-4)t^3
\end{split}
\end{align}
such that $\dv(f''/x^3)=C''+{}^\tau\!C''+{}^{\tau\tau}\!C''-3H$.
Continuing, we obtain a rational function whose norm times
$(-15\zeta^2/2)f''/x^3$ is equal to $f/x^3$.

Therefore the rational functions
\begin{equation}
\label{eqn:multiplerepresentatives}
\frac{f}{x^3},\qquad
\frac{\zeta}{4}\,\frac{f'}{x^3},\qquad
\frac{-15\zeta^2}{2}\,\frac{f''}{x^3}
\end{equation}
have the property that their corresponding Azumaya algebras
(for the cyclic extension
$K_0$ of $k$) all define (restrictions of) the same element
of $\Br(X)$.
Each Azumaya algebra is defined over an open subset of $X$.
We would like to have $X$ covered by such open sets.
Indeed, the divisors on $X$ defined by
$f$, $f'$, and $f''$ have trivial intersection.
All we do now is replace $H$ by other
hyperplane sections $H': y=0$, etc.
So, in \eqref{eqn:multiplerepresentatives} we replace
$x$ by the other coordinate functions, and obtain a larger collection
of Azumaya algebras which all represent the same element of $\Br(X)$ and
whose domains of definition cover $X$.
\end{exam}

\section{Local invariants}
\label{sec:localinv}

The classical Nullstellensatz theorem asserts that if
polynomials $f_1$, $\ldots$, $f_m\in k[x_1,\ldots,x_N]$ define
the empty scheme in affine space $\A^N_k$ then there exist polynomials
$g_i$ satisfying $\sum_i f_ig_i=1$.
An effective arithmetic Nullstellensatz theorem \cite{kps}
applies to a number field $k$ with ring of integers $\mathfrak{o}_k$.
Then, assuming $f_i\in \mathfrak{o}_k[x_1,\ldots,x_N]$
for all $i$, there exist
$g_i\in \mathfrak{o}_k[x_1,\ldots,x_N]$
satisfying
\begin{equation}
\label{eqn:nss}
\sum_i f_ig_i=\varpi
\end{equation}
for some $\varpi\in \mathfrak{o}_k$.
The theorem asserts that such $g_i$ can be found, satisfying
bounds on their degrees and on
the heights of their coefficients.
The bounds are effective and depend on
the degrees and the heights of the coefficients of the $f_i$.

The given Brauer group element $A$ is unramified at all but a
finite set of places of $k$.
Consider a fixed integral model $\mathcal{X}$ of $X$ over
$\mathfrak{o}_k$.
Recall, a place $v$ of $k$ is said to be
of good reduction when the integral model is
smooth over the residue field ${\mathbf k}_v$.
At all but finitely many of these places, the polynomials $f_t$ and $g_t$
appearing in the cocycles reduce to nontrivial elements
of $\cO(\mathcal{X}_{{\mathbf k}_v})$.
Then by purity \cite[(III.6.1)]{gb}, the Brauer group elements are
unramified at all $v$-adic points of $X$.

Consequently, it is only necessary to carry out the local analysis,
i.e., the computation of local invariants \eqref{eq:localinvariants}
of $X(k_v)$-points, at those $v$ which are
places of bad reduction, places where the
rational functions appearing in the cocycles fail to extend,
and real places of $k$.

We first treat a non-archimedean place $v$.
Choose a valuation $v_1\mid v$ of $K$ and an embedding
$K\to K_{v_1}$.
The Galois group of $\Gal(K_{v_1}/K)$ is the subgroup $G_1$ of
elements of $G$ which stabilize $v_1$.
The embedding $K^{G_1}\to K_{v_1}$ then factors through
$k_v$.
Replacing $k$ by $K^{G_1}$ we are thereby reduced to the
case that $v$ extends uniquely to a valuation $v_1$ of $K$, and
$\Gal(K_{v_1}/k_v)=G$.
We then have $G$-equivariant inclusions
$$\cO(U_K)\to \cO(U_{K_{v_1}})$$
for any variety $U$ over $k$.

We fix an element $A\in \Br(X)$, an open covering $\{U^j\}$ of $X$,
and $2$-cocycle representatives
$\{\beta^j\}$ with values in $\cO(U^j_K)^*$, as constructed in
Steps 1 through 5
(so that $\beta^j$ is a representative of the
restriction of $A$ to $U^j$ for each $j$).
We may suppose each $\beta^j$ is a tuple $(f^j_t/g^j_t)_{t\in T}$
for some finite index set $T$,
indexing a basis of $\mathsf{C}^2(M)$.
We may require that $f^j_t$ and $g^j_t$ should be polynomials with
coefficients in the ring of integers $\mathfrak{o}_K$.
These are constructed algorithmically, so it is possible to give an
effective bound on the valuation of the coefficients.

\begin{lemm}
\label{nbhd1trivial}
Let $k$ be a number field, $K$ a finite Galois extension of $k$ with
Galois group $G$, and let $v$ be a non-archimedean valuation of $k$
admitting a unique extension to a valuation $v_1$ on $K$.
Fix a resolution $\mathsf{C}^\bullet(K_{v_1}^*)$ for computing group cohomology of
the $G$-module $K_{v_1}^*$, and let us write elements of
$\mathsf{C}^2(K_{v_1}^*)$ as tuples $(\alpha_t)$ indexed by $t\in T$ for some
finite index set $T$.
Then there is, effectively in terms of the given data,
a number $N$ such that
every $2$-cocycle $(\alpha_t)$ with
$\min_t v_1(\alpha_t-1)>N$ is a $2$-coboundary.
\end{lemm}

\begin{proof}
By restriction of scalars, there is a smooth affine variety $Z$ over $k_v$
such that $2$-cocycles $(\alpha_t)$ with values in $K_{v_1}^*$
map bijectively to $Z(k_v)$,
where the map is $k_v$-linear and explicit.
Let $z_1\in Z$ correspond to the $2$-cocycle $(1)$.
There is, furthermore, a smooth affine variety $Y$ over $k_v$ and
a similar bijection with $1$-cochains, and a smooth surjective map of
$k_v$-varieties $\varrho\colon Y\to Z$ such that $2$-coboundaries
correspond to points of $\varrho(Y(k_v))$.

Fix coordinates for $Y$ and $Z$.
Equations for $Y$, $Z$, and $\varrho$ can be written down explicitly.
By the Inverse Function Theorem, there is an explicit
$v$-adic neighborhood $W$ of $z_1$ in $Z(k_v)$ and an analytic splitting
(not needed explicitly) of $\varrho\colon Y(k_v)\to Z(k_v)$ on $W$.
Hence for every $k_v$-point of $W$ the corresponding $2$-cocycle $(\alpha_t)$
is a $2$-coboundary.

The restriction of scalars map, sending a $K_{v_1}^*$-valued cocycle
$(\alpha_t)$ to a point $z\in Z(k_v)$, is given coordinatewise by
$k_v$-linear expressions.
Since $v_1$ is Galois-invariant, each coordinate of $z$ has
$v$-adic valuation bounded from below by
$\min_t(v_1(\alpha_t))$ plus an explicit constant.
So there is an effective bound $N$ such that
$\min_t(v_1(\alpha_t-1))>N$ implies $z\in W$, so that $(\alpha_t)$ must be
a $2$-coboundary.
\end{proof}

\begin{coro}
\label{cor:closecocycles}
Keep the notation of Lemma \ref{nbhd1trivial}.
Fix an integer $P$, and
let $(\alpha_t)$ and $(\alpha'_t)$
be $2$-cocycles with values in $K_{v_1}^*$ with
$|v_1(\alpha_t)|\le P$ and $|v_1(\alpha'_t)|\le P$ for all $t$.
Then there is an effective bound $Q$, depending on the fields
$k_v$ and $K_{v_1}$ and $P$ but not on the given $2$-cocycles,
such that
if $v_1(\alpha_t-\alpha'_t)\ge Q$ for all $t$
then $(\alpha_t)$ and $(\alpha'_t)$ define the same element of
$\Br(k_v)$.
\end{coro}

The proof of the next result contains an algorithmic description of
the computation of local invariant in $\Q/\Z$ of the element in
$\Br(k_v)$, given as a $2$-cocycle $(\alpha_t)$.

\begin{prop}
\label{computelocalinvariant}
With notation as in Lemma \ref{nbhd1trivial}, let an integer $P$ be
given, and let
$(\alpha_t)$ be a $2$-cocycle with values in $K_{v_1}^*$ such that
$|v_1(\alpha_t)|\le P$ for all $t$.
Then there is an effective computation, taking as input the collection of
$\alpha_t$ each specified to an effectively determined degree of precision,
of the local invariant of the element of $\Br(k_v)$ corresponding to this
$2$-cocycle.
\end{prop}

\begin{rema}
\label{itimpliesthelemma}
Formally speaking, Proposition \ref{computelocalinvariant} implies
Lemma \ref{nbhd1trivial}.
We have included Lemma \ref{nbhd1trivial} because it admits a direct
proof independent of the detailed algorithmic treatment of the proof
of Proposition \ref{computelocalinvariant}.
Lemma \ref{nbhd1trivial} is significant because it can be
used to simplify the computation in practice.
Indeed
the bound $N$ arising in the proof of the lemma can be quite reasonable,
while the algorithmic description tends to lead to a much worse bound.
See the example at the end of this section, where such bounds are
obtained, compared, and used to help complete the local analysis.
\end{rema}

\begin{proof}[Proof of Proposition \ref{computelocalinvariant}]
It is possible to obtain an explicit map between the
given resolution $\mathsf{C}^\bullet$ and the standard resolution
$\mathsf{C}^\bullet_G$ with
$$\mathsf{C}^2_G(M)=\bigoplus_{(g,h)\in G\times G} M.$$
Thus we are reduced to the case of given cocycle data
$(\alpha_{g,h})$ satisfying the standard cocycle condition
$$\alpha_{g,h}\alpha_{gh,j}=\alpha_{g,hj}{}^g\alpha_{h,j}$$
for all $g$, $h$, $j\in G$.
There is no loss of generality in assuming that for all $g\in G$ that
$\alpha_{e,g}=\alpha_{g,e}=1$, where $e$ denotes the identity element of $G$.
Then the identification via crossed products
$$H^2(\Gal(\bar k_v),\bar k_v^*)=\Br(k_v)$$
associates to $(\alpha_{g,h})$ the element of the Brauer group
corresponding to the central simple algebra
$$\cA=\bigoplus_{g\in G} K_{v_1}\cdot e_g$$
where $e_g$ are generators, subject to the relations
\begin{gather*}
e_g e_h = \alpha_{g,h} e_{gh}\\
e_g c = {}^g\!c e_g
\end{gather*}
for $g$, $h\in G$ and $c\in K_{v_1}$ (see, e.g., \cite{draxl}).

Now we recall the local invariant isomorphism $\Br(k_v)=\Q/\Z$.
There exists, up to isomorphism, a unique unramified extension
$\ell_v/k_v$ of each degree $d$;
such $\ell_v$ is cyclic Galois over $k_v$.
Taking $d$ such that the class of $\cA$ in $\Br(k_v)$ is $d$-torsion
(e.g., we can take $d=|G|$), there must exist an isomorphism
\begin{equation}
\label{eq:chiiso}
\chi\colon \cA\otimes_{k_v} \ell_v\to M_{|G|}(\ell_v)
\end{equation}
of the extension of $\cA$ to the degree $d$ unramified field extension $\ell_v$
with the $|G|$-by-$|G|$ matrix algebra over $\ell_v$.
The extension $\ell_v$ can be obtained explicitly
(e.g., from a cyclotomic extension).

We let $n=|G|$, and write
$\cA_{k'}$
for $\cA\otimes_{k_v} k'$ where $k'/k_v$ is any field extension.
The first step toward the construction of
an isomorphism \eqref{eq:chiiso} is to
find an $\ell_v$-point on the Brauer--Severi variety
associated with $\cA$.
The Brauer--Severi variety of $\cA$ has the property that
$k'$-points correspond bijectively with
$n$-dimensional left ideals of $\cA_{k'}$, for any extension $k'/k_v$.
This description supplies explicit equations for
the Brauer--Severi variety as a subvariety of the Grassmannian variety
$Gr(n,n^2)$ over $k_v$ (see, e.g., \cite{artin}), and hence there are
\emph{a priori} bounds on the $v$-adic height of an $\ell_v$-valued solution.
We let
$\cA_1\subset \cA_{\ell_v}$ denote the corresponding left ideal.

There are $1=x_1$, $x_2$, $\ldots$, $x_n\in \cA_{\ell_v}$ such that
if we set $\cA_i=\cA_1x_i$ then the spaces $\cA_1$, $\ldots$, $\cA_n$ span $\cA_{\ell_v}$.
The $\cA_i$ can be found effectively, and their entries can be effectively
bounded.
More precisely, if we vary over a choice of
possible $x_i$ (e.g., a fixed
$k_v$-basis of $\cA$), then for each choice we get a map
$(\cA_1)^n\to \cA$, given by $(y_i)\mapsto \sum_i y_ix_i$.
Such a map is represented by a matrix; now we consider the determinant
of the matrix.
The tuple of determinants over all choices of $x_i$ is nonvanishing,
so appealing to an effective Nullstellensatz we have an \emph{a priori} bound
on the minimum of the determinants.

Under the isomorphism of vector spaces
$$\cA_{\ell_v}=\cA_1\oplus \cdots \oplus \cA_n$$
the vector $1\in \cA_{\ell_v}$ maps to some
$(e_{1,1},\ldots,e_{n,n})\in \cA_1\oplus\cdots\oplus \cA_n$.
Now we define
$\chi$ to map $e_{i,i}$ to the matrix with a single entry $1$ in the
$i$th row and $i$th column, and all other entries $0$.
The prescription
$e_{i,j}\in e_{i,i}\cA_{\ell_v}\cap \cA_{\ell_v}e_{j,j}$,
defines the $e_{i,j}$ uniquely up to scale, and we can thereby
choose $e_{i,i+1}$ for each $i$.
Then the relations
$e_{i,j}e_{i',j'}=\delta_{j,i'}e_{i,j'}$ determine all the $e_{i,j}$
and complete the definition of the algebra isomorphism $\chi$.

The action of $g\in G$ on $\ell_v$ induces the algebra automorphism
$\varphi_g\colon M_n(\ell_v)\to M_n(\ell_v)$ defined by
$$\varphi_g(m)={}^g(\chi({}^{g^{-1}}(\chi^{-1}(m))))$$
Write
$$\varphi_g(m)=p_g^{-1}mp_g.$$
This is effective Skolem--Noether:
for any vector $w$, the matrix $q$ with vector $\varphi(e_{i,1})w$ in
column $i$ (where $e_{i,1}$ denotes the matrix with $1$ in row $i$ column $1$
and $0$ elsewhere) for each $i$ satisfies
$\varphi(m)q=qm$ for all $m\in M_n(\ell_v)$, and
for some choice of $w$ the determinant of $q$ must
have valuation less than some bound that can be made explicit.
Then take $p_g=q^{-1}$.

Now we follow \cite[\S X.5]{serre}: if we define $\beta_{g,h}$ by
\begin{equation}
\label{betagh}
\beta_{g,h}^{-1}=p_g {}^g\!p_h p_{gh}^{-1}
\end{equation}
then $(\beta_{g,h})$ is a $2$-cocycle with values
in $\ell_v^*$ for the extension $\ell_v/k_v$
such that its class in $\Br(k_v)$ is equal to $\cA$.
(As remarked in \emph{ibid}., the construction applied to an algebra defined
by a $2$-cocycle for the extension $\ell_v$ representing
$B\in H^2(\Gal(\bar k_v/k_v),\bar k_v^*)=\Br(k_v)$
naturally gives rise to a $2$-cocycle for the class $-B$; this explains the
power $-1$ in the definition of $\beta_{g,h}$.)
Now $v(\beta_{g,h})$ is a $2$-cocycle with values in $\Z$.
Using the long exact sequence of cohomology of the sequence
$$0\to \Z\to \Q\to \Q/\Z$$
and acyclicity of $\Q$ we obtain a $1$-cocycle with values in $\Q/\Z$,
that is, a group homomorphism
$$\Gal(\ell_v/k_v)\to \Q/\Z.$$
The local invariant of the algebra
$\cA$ (associated with the given $2$-cocycle
with values in $K_{v_1}^*$)
is the image in $\Q/\Z$ of the unique element of
$\Gal(\ell_v/k_v)$ which induces the Frobenius automorphism on
residue fields (cf.~\cite[\S XIII.3]{serre}).
\end{proof}

We write $k_v$-points of $X$ as tuples
$$\mathbf{x}=(x_0,\ldots,x_N)$$
with $x_{i_0}=1$ for some ${i_0}$ and $x_i\in \mathfrak{o}_v$
(the ring of integers of $k_v$) for all $i$.
It suffices to consider one chart at a time (i.e., choice of $i_0$).
Henceforth we work with points in a single chart.

Set $h^j:=\prod_t f^j_tg^j_t$.
The $h^j$ are nonvanishing on $U^j$, which cover $X$.
So we can apply
effective Nullstellensatz \eqref{eqn:nss} to the $h^j$ to deduce that
there are homogeneous polynomials $H^j$ with
coefficients in $\mathfrak{o}_v$ and
an equation
$$\sum_j h^j H^j = A x_{i_0}^r$$
for some $A\in \mathfrak{o}_v$ and $r$.
The valuations of $A$ and of the coefficients of the $H^j$
are effectively bounded, as are the degrees of the $H^j$.
It follows that for each $\mathbf{x}\in X(k_v)$
there exists $j$
such that the cocycle $(\beta^j_t)$ can be evaluated at $\mathbf{x}$,
with $\max(v(f^j_t(\mathbf{x})),v(g^j_t(\mathbf{x})))\le P$
for some effective bound $P$ that is independent of $\mathbf{x}$.

\begin{prop}
\label{closepointsofX}
Fix the data of Assumption \ref{assumption1}
and consider $A\in \Br(X)$ with representing cocycles
$(\beta^j_t)$ from Steps 1 through 5.
There is an effective bound $Q'$ such that if
$\mathbf{x}$ and $\mathbf{x}'$ are points of
(a fixed chart of) $X(k_v)$ satisfying
$v(x_i-x'_i)>Q'$ for all $i$,
then $\inv_v(A(\mathbf{x}))=\inv_v(A(\mathbf{x}'))$.
\end{prop}

\begin{proof}
There is an effective bound $P$ such that
for some $j$ we have
$$\max(v(f^j_t(\mathbf{x})),v(g^j_t(\mathbf{x})))\le P$$
for all $t$.
Hence, for large enough $Q'$,
we also have
$\max(v(f^j_t(\mathbf{x}')),v(g^j_t(\mathbf{x}')))\le P$
for this value of $j$ and all $t$.
The result now follows by Corollary \ref{cor:closecocycles}.
\end{proof}

We can now give
an effective procedure to compute the local invariants of
the $k_v$-points of $X$.
As before, we focus on a single chart of $X$.
Then, with $Q'$ as in
Proposition \ref{closepointsofX},
we enumerate the classes of points of
$X(k_v)$, where points $\mathbf{x}$, $\mathbf{x}'$
are considered to lie in the same class if we have
$v(x_i-x'_i)>Q'$.
For each class, we choose a representative point and compute it to
a precision which is
\begin{itemize}
\item[(i)] enough to identify the $j$ such that
$\max(v(f^j_t(\mathbf{x})),v(g^j_t(\mathbf{x})))\le P$ for all $t$,
\item[(ii)] enough to carry out the algorithm of
Proposition \ref{computelocalinvariant}.
\end{itemize}
Then we evaluate $(\beta^j_t(\mathbf{x}))$ and apply the algorithm of
Proposition \ref{computelocalinvariant}.
The output is the local invariant of all points
$\mathbf{x'}\in X(k_v)$ lying in the same class (in the sense
we have just introduced) as $\mathbf{x}$.

\medskip

We now treat the case of a real place
$v$ of $k$.
Then $X(k_v)$ is a real algebraic manifold, and the
local invariant of the Brauer group element $A$ is
constant on connected components.
So we are reduced to determining whether none, some, or all of the
connected components of $X(k_v)$ have points where a
cocycle representative of $A$ has values in the image of an
algebraic coboundary map.
For real algebraic varieties there are effective procedures
for computing the number of connected components
and the image of an algebraic map
(cf.\ \cite{hrr}).
So the ramification pattern of any Brauer group element
(or collection of Brauer group elements) can be determined effectively.

\begin{exam}
We carry out the local analysis in the Cassels--Guy example
\begin{equation}
\label{eq:591012}
5x^3+9y^3+10z^3+12t^3=0
\end{equation}
to deduce that this cubic surface $X$ violates the Hasse principle.
The original proof of this fact from \cite{cg} is based on
ideal class group computations.
The work of Colliot-Th\'el\`ene, Kanevsky, and Sansuc \cite{ctks} treated
diagonal cubic surface using the Brauer--Manin obstruction; however they
used more complicated (non-cyclic) Azumaya algebras.

The base field is $k=\Q(\zeta)$ with $\zeta=e^{2\pi i/3}$.
We have produced an element $A\in \Br(X)$ which generates
$\Br(X)/\Br(k)\cong \Z/3\Z$.
On various open subsets, $A$ is represented by cyclic Azumaya algebras,
for the field extension $k(\sqrt[3]{2/3})$, of the function field elements
$$g_1:=\frac{f}{x^3},\qquad
g_2:=2\zeta\frac{f'}{x^3},\qquad
g_3:=-60\zeta^2\frac{f''}{x^3}$$
where $f$, $f'$, and $f''$ are the polynomials of
\eqref{eq:vanishonCtCttC},
\eqref{fprime}, and \eqref{fprimeprime}.

We proceed with the evaluation of the functions $g_i$ at $v$-adic points of
$X$.
The equation \eqref{eq:591012} has good reduction outside of the primes $2$,
$\sqrt{-3}$, and $5$.
Since the functions $g_i$ reduce to nontrivial rational functions
on the reduction at any place $v$ outside these primes,
it follows that $A$ is unramified at any point of $X(k_v)$.

Consider first $v=2$.
We proceed to evaluate the functions $g_i$ at $2$-adic points of $X$
and compute the local invariants of the resulting $2$-cocycles.
Every point of $X(k_2)$ with $2$-adic integer coefficients must
have $x$ and $y$ not divisible by $2$, so we are free to consider
$x=1$ throughout the analysis.
It suffices to consider $y$, $z$, and $t$ mod $8$, as we see in the
following computations.
For all $y$, $z$, and $t$
satisfying \eqref{eq:591012} mod $8$,
there is always some $i$ such that $g_i$
is of the form
$2^n(1+2a)$ with $a\in \mathfrak{o}_2$ and $0\le n\le 2$.
Any number of this form is a norm, since $2$ is the norm of
$2+\sqrt[3]{12}+\sqrt[3]{18}$
and for any $a\in \mathfrak{o}_2$ the number $1+2a$ is a cube in $k_2$.
So the corresponding $2$-cocycle is always trivial, and the local invariant
at all $2$-adic points of $X$ is $0$.

Now consider the place $v=\sqrt{-3}$.
Again, we may suppose $x=1$.
We observe that at all points of $X(k_v)$ the value of $g_1$ is
$\sqrt{-3}$ times a unit $u\in \mathfrak{o}^*_v$.
To determine the required precision, we make some remarks concerning
norms for the extension $k_v(\sqrt[3]{2/3})/k_v$.
We have that $\sqrt{-3}$ is a norm,
$$\sqrt{-3}=N\bigl(-1-2\zeta+(1-\zeta)\sqrt[3]{2/3}\bigr),$$
and that any $u\in \mathfrak{o}^*_v$ congruent to $1$ mod $9$
is a cube in $k_v$.
Consequently it suffices to evaluate $(1/\sqrt{-3})g_1$ modulo $9$,
and for this it suffices to consider $y$, $z$, and $t$ modulo $9\sqrt{-3}$.
It results that $(1/\sqrt{-3})g_1$, modulo $9$, always equals
one of the following:
\begin{equation}
\label{eq:6choices}
\zeta,\quad 4\zeta,\quad 7\zeta,\quad 3+\zeta, \quad 3+4\zeta,\quad
3+7\zeta.
\end{equation}
The equation from the previous paragraph expressing $2$ as a norm is
valid in $k_v$, and we also have
$$N\bigl(1+(-1-2\zeta\bigr)\sqrt[3]{2/3})=3+4\zeta.$$
Consequently, all numbers in \eqref{eq:6choices} are of the form
$\zeta$ times a norm.
By the algorithm of Proposition \ref{computelocalinvariant},
the local invariant of the cyclic division algebra defined by $\zeta$
is computed to be $2/3$.
So the local invariant
at all $v$-adic points of $X$ is $2/3$.

We remark that the application of Proposition \ref{computelocalinvariant}
in the previous paragraph to compute a $v$-adic local invariant
requires solving for
the isomorphism \eqref{eq:chiiso} to certain precision in order to
ensure that the subsequent steps of the algorithm
lead to a meaningful evaluation.
The observed loss of significance in computing $\beta_{g,h}$
\eqref{betagh} from $\chi$
is on the order of $3^{10}$, meaning
that according to the algorithmic description of the computation
we would require $y$, $z$, and $t$ to such a high precision
(the \emph{a priori} bounds are of course much worse).
This is in contrast to the analysis from Lemma \ref{nbhd1trivial},
which implies that values modulo $9\sqrt{-3}$ suffice.

Finally,
the field $k_5$ contains $\sqrt[3]{2/3}$, so the $5$-adic analysis is trivial.
At any adelic point of $X$, the sum of local invariants is $2/3$.
In conclusion,
the surface $X$ defined by \eqref{eq:591012} violates the Hasse principle.
\end{exam}


\end{document}